\documentclass[11pt, a4paper,leqno]{amsart}
\usepackage{amsmath,amsthm,amscd,amssymb,amsfonts, amsbsy}
\usepackage{latexsym}
\usepackage{txfonts}
\usepackage{exscale}
\usepackage{mathrsfs}
\usepackage[titletoc]{appendix}

\usepackage[utf8]{inputenc}

\renewcommand\labelenumi{(\roman{enumi})}
\renewcommand\theenumi\labelenumi

\usepackage[colorlinks,linkcolor=blue,citecolor=blue,pagebackref,hypertexnames=false, breaklinks]{hyperref}
\usepackage[usenames,dvipsnames]{color}

\calclayout
\allowdisplaybreaks

\def \Int{\,\Int\,}

    \def\XXint#1#2#3{{\setbox0=\hbox{$#1{#2#3}{\int}$}
    \vcenter{\hbox{$#2#3$}}\kern-.5\wd0}}

    \def\XXint#1#2#3{{\setbox0=\hbox{$#1{#2#3}{\int}$}
    \vcenter{\hbox{$#2#3$}}\kern-.5\wd0}}

\def \C{ \mathbb{C} }

\def \R{ \mathbb{R} }

\renewcommand{\chi}{\mathbf{1}}

\theoremstyle{plain}
\newtheorem{theorem}[equation]{Theorem}
\newtheorem{lemma}[equation]{Lemma}
\newtheorem{corollary}[equation]{Corollary}

\theoremstyle{definition}
\newtheorem{definition}[equation]{Definition}

\newtheorem{example}[equation]{Example}

\theoremstyle{remark}
\newtheorem{remark}[equation]{Remark}

\numberwithin{equation}{section}

\numberwithin{equation}{section}

\numberwithin{equation}{section}

\begin{document}
\allowdisplaybreaks
\author{Alberto Lastra}
\address{Alberto Lastra
\\
Universidad de Alcalá 
\\
Departamento de Física y Matemáticas
\\
Campus universitario 
\\E-28805 Alcalá de Henares (Madrid), Spain
}
 \email{alberto.lastra@uah.es}
\author{Cruz Prisuelos-Arribas}
\address{Cruz Prisuelos-Arribas
\\
Universidad de Alcalá 
\\
Departamento de Física y Matemáticas
\\
Campus universitario 
\\E-28805 Alcalá de Henares (Madrid), Spain
} \email{cruz.prisuelos@uah.es}

\title[Linear systems of moment differential equations via generalized matrix exponentials]{Solutions of linear systems of moment differential equations via generalized matrix exponentials}

	\thanks{The work of the first author is partially supported by the project PID2019-105621GB-I00 of Ministerio de Ciencia e Innovación, Spain, and by Dirección General de Investigación e Innovación, Consejería de Educación e Investigación of the Comunidad de Madrid (Spain), by Universidad de Alcalá under grant CM/JIN/2021-014, Proyectos de I+D para Jóvenes Investigadores de la Universidad de Alcalá 2021, and the Ministerio de Ciencia e Innovación-Agencia Estatal de Investigación MCIN/AEI/10.13039/501100011033 and the European Union ``NextGenerationEU''/ PRTR, under grant TED2021-129813A-I00.
	The work of the second author is supported by the project	
	PID2019-­107914GB-­I00 of Ministerio de Ciencia e Innovación, Spain.}


\date{\today}

\begin{abstract}
A generalized exponential matrix based on the construction of kernel operators for generalized summability is defined and analyzing its main properties, generalizing the classical exponential matrix and fractional exponential matrix. 

This object serves as a practical tool to express the solutions of linear systems of moment differential equations in a compact manner, in the spirit of the classical exponential matrix.
\end{abstract}

\maketitle

\tableofcontents

\section{Introduction}\label{secintro}
In this paper we consider the linear system of moment differential equations of the form
\begin{align}\label{problem}
\partial_my(z)=Ay(z),
\end{align}
where $A\in \mathbb{C}^{n\times n}$ is a constant matrix and $y(z)=(y_1(z),\ldots,y_n(z))$ is a vector of unknown functions, for some positive integer $n\geq 1$ and $\partial_my(z)=(\partial_m(y_1),\ldots,\partial_m(y_n))$ denotes the moment derivative of $y(z)$ (see Section~\ref{secmomder}). The explicit general solution to (\ref{problem}) was stated in~\cite{Lastra22}. The main aim of the present work is to provide a practical tool to represent the solutions of such system. This will be done following the spirit of the exponential matrix in the classical setting, i.e. when a linear system of first order differential equations is considered. Our main aim is to express the general solution of (\ref{problem}) in terms of a generalized exponential matrix $E(Az)$. 

The importance of moment differential equations has increased in recent years due to the versatility of the moment differential operator which generalizes the classical derivative, but also Caputo fractional derivative and also $q-$derivative, when choosing the appropriate moment sequence. The first step in this theory is due to W. Balser and M. Yoshino~\cite{bayo} stating the formal definition of moment derivation. In the last decade, this definition has been extended not only to holomorphic functions on a neighborhood of the point, but also to sums and multisums of formal power series (see~\cite{lamisu1,lamisu2}) which also show recent advances achieved on the summability of the formal solutions to such functional equations, and also on the estimation of upper bounds for the coefficients of their formal solutions (see~\cite{lamisu0,sa,su}) or the Stokes phenomenon (see~\cite{mitk}).

It is known that the choice of moment sequence $m_1=(p!)_{p\ge0}$ provides usual derivation, whereas fixing $s>0$, the moment derivative associated to the sequence $m_s=(\Gamma(1+\frac{p}{s})_{p\ge0}$ is quite related to Caputo fractional derivative. More precisely, $\partial_{m_s}y(z^{1/s})={}^{C}D_z^{1/s}(y(z^{1/s}))$, where $^{C}D_z^{1/s}$ stands for Caputo fractional derivative of order $1/s$. For every fixed $q>1$, the sequence $m_q=(\Gamma_q(1+p))_{p\ge0}$ determines $q-$derivatives, where $\Gamma_q(\cdot)$ is the $q$-analog of Gamma function (see~\cite{zh}, as a reference in this sense). As a matter of fact, $m$ is said a moment sequence due to $m_s$ for $s>0$ turns out to be the sequence of moments associated to a positive Laguerre-type measure, and $m_q$ is intimately related to the sequence $(q^{p(p-1)/2})_{p\ge0}$ in the sense that their classical ultradifferentiable classes of functions coincide, and this last sequence is the sequence of moments associated to a measure of usual appearance in summability theory (see~\cite{drlama,lama,razh} as examples of this fact). Nevertheless, the construction of the generalized exponential matrix (see~\ref{exponential}) can be generalized to any sequence of positive real numbers $m$ in a merely formal way, but also providing an entire function under a mild assumption on the sequence $m$ (see (\ref{e433})). We will distinguish these important cases for applications in the paper.

The construction of matrix exponentials is of great importance in applications. The classical construction of the exponential of a matrix is a widespread tool to solving linear systems of differential equations. In~\cite{ro}, the author introduces the concept of fractional matrix exponential which solves linear systems of fractional differential equations of the form
$$D^{\alpha}y(z)=A^{\alpha}y(z),$$
where $D^{\alpha}$ is a fractional differential operator, and $\alpha>0$. The matrix $A^{\alpha}$ generalizes natural powers of the matrix $A$. Indeed, the solution is formulated in the form $y(z)=\exp(Az)$, constructed in two cases: for $D^{\alpha}$ being Caputo fractional derivative or Riemann-Liouville fractional differential operator. 

\vspace{0.3cm}

\textbf{Notation:}

For all $d\in\R$ and $0<\theta<2\pi$, $S_{d}(\theta)$ stands for the complex numbers $z$ such that $|\arg(z)-d|<\frac{\theta}{2}$.

For all $z_0\in\C$ and $r>0$, $D(z,r)$ is the open disc centered at $z_0$ and with radius $r$.

Let $\mathbb{E}$ be a complex Banach space. We denote $\mathbb{E}[[z]]$ the space of formal power series with coefficients in $\mathbb{E}$. $\mathbb{E}\{z\}$ stands for the set of power series with positive radius of convergence on some neighborhood of the origin.

\section{Preliminaries and main definitions}

In a first subsection, we provide the main elements present in the constructions of kernels for generalized summability, and moment differentiation. A second subsection is devoted to state the definition of the generalized moment exponential matrix, followed by its main properties. The generalized exponential of Jordan block decomposition of a matrix provides the final structure of the generalized moment exponential matrix in the last two subsections. We have split the discussion into diagonalizable and non-diagonalizable matrices.

\subsection{Strongly regular sequences and generalized summability}\label{secmomder}

The definitions and results in this section can be found in detail in the seminal work by J. Sanz~\cite{sa0} on generalized summability, and the references therein. 

Following V. Thilliez~\cite{th}, a strongly regular sequence is defined as follows.

\begin{definition}
A sequence of positive real numbers $\mathbb{M}=(M_p)_{p\ge0}$ with normalization $M_0=1$ is said to be a strongly regular sequence if it satisfies:
\begin{itemize}
\item[(lc)] $M_p^2\le M_{p-1}M_{p+1}$, for all $p\ge1$ (i.e. $\mathbb{M}$ is a logarithmically convex sequence).
\item[(mg)] There exists $A>0$ such that $M_{p+q}\le AM_pM_q$ for all $p,q\ge0$ (i.e. $\mathbb{M}$ is of moderate growth).
\item[(snq)] There exists $A>0$ such that 
$$\sum_{q\ge p}\frac{M_q}{(q+1)M_{q+1}}\le A\frac{M_p}{M_{p+1}},\quad p\ge0,$$
(i.e. $\mathbb{M}$ is non-quasianalytic).
\end{itemize}
\end{definition}

The most outstanding example of strongly regular sequence is $\mathbb{M}_{\alpha}=(p!^{\alpha})_{p\ge0}$, for a fixed $\alpha>0$, known as Gevrey sequence of order $\alpha$.

Given a strongly regular sequence $\mathbb{M}$, the function $M:[0,\infty)\to [0,\infty)$ defined by $M(0)=0$ and $M(t)=\sup_{p\ge0}\log\left(\frac{t^p}{M_p}\right)$ for $t>0$, measures the growth of $\mathbb{M}$ and allows to define the positive number $\rho(M)=\limsup_{r\to\infty}\frac{\log(M(r))}{\log(r)}$.

\begin{definition}
Let $\mathbb{M}$ be a strongly regular sequence, and let $M$ and $\rho(M)$ be defined from $\mathbb{M}$ as before. The sequence $\mathbb{M}$ admits a pair of kernel functions for $\mathbb{M}$-summability if there exist $e$ and $E$ such that:
\begin{enumerate}
\item $e\in\mathcal{O}(S_0(\frac{\pi}{\rho(M)}))$. Moreover, there exists $t_0>0$ and for all $z_0\in S_{0}(\frac{\pi}{\rho(M)})$ there exists $r_0=r_0(z_0)>0$ with $D(z_0,r_0)\subseteq S_{0}(\frac{\pi}{\rho(M)})$ such that 
$$\int_0^{t_0}\sup_{z\in D(z_0,r_0)}|e(t/z)|\frac{dt}{t}l\infty.$$
We also assume that for every $\delta>0$ there exist $K_1,K_2>0$ with
$$|e(z)|\le K_1\exp\left(-M(|z|/K_2)\right),\quad z\in S_0(\frac{\pi}{\rho(M)}-\delta).$$
\item $e(x)$ is a positive real number for all $x\in\R$ with $x>0$.
\item $E\in\mathbb{O}(\C)$, and
$$|E(z)|\le K_3\exp\left(M(|z|/K_4)\right),\quad z\in\C,$$
for some $K_3,K_4>0$. Additionally, there exists $\beta>0$ such that for all $0<\theta<2\pi-\frac{\pi}{\rho(M)}$ and $R>0$, there exists $K_5>0$ with $|E(z)|\le K_5/|z|^{\beta}$ for $z\in S_{\pi}(\theta)\cup(\C\setminus D(0,R))$.
\item Taylor expansion of $E$ at the origin is given by
$$E(z)=\sum_{p\ge0}\frac{z^p}{m(p)},\quad z\in\C,$$
where $m=(m(p))_{p\ge0}$ is the sequence of Laguerre-type moments of the function $e$, i.e.
\begin{equation}\label{e336}
m(p)=\int_{0}^{\infty}t^{p-1}e(t)dt,\qquad p\ge0.
\end{equation}
\end{enumerate}
\end{definition}

Given a strongly regular sequence $\mathbb{M}$, a sufficient condition for the existence of an associated pair of kernel functions for generalized summability is that $\mathbb{M}$ admits a nonzero proximate order (see~\cite{jisasc,lamasa}), which is the case of most sequences appearing in applications, in particular for Gevrey sequence $\mathbb{M}_{\alpha}$, for every $\alpha>0$. 

Let $m=(m_p)_{p\ge0}$ be a sequence of positive real numbers (we will usually consider $m$ to be the sequence of moments associated to a the kernel functions for generalized summability associated to a given strongly regular sequence). In 2010, W. Balser and M. Yoshino~\cite{bayo} put forward the concept of moment differentiation as the operator $\partial_{m,z}:\C[[z]]\to\C[[z]]$ defined by
$$\partial_{m,z}\left(\sum_{p\ge0}\frac{a_p}{m_p}z^p\right)=\sum_{p\ge0}\frac{a_{p+1}}{m_p}z^p.$$

We will usually denote the moment differential operator by $\partial_{m}$, as no confusion with the variable may arise in this work.

This definition can be naturally extended to formal power series with coefficients in a complex vector space. In the present study, we will usually consider the coefficients in the Banach space of $n\times n$ matrices of complex entries, $\C^{n\times n}$, endowed with a submultiplicative norm. This definition makes sense when considering holomorphic functions defined on some neighborhood of the origin, by identifying the function with its Taylor expansion at the origin. The most outstanding situation is that of $\mathbb{M}_{\alpha}=(p!^\alpha)_{p\ge0}$. In this case, the sequence of moments is given by $m_{k}=(\Gamma(1+\frac{p}{k})_{p\ge0}$ and the moment derivative is directly linked to Caputo fractional derivative: $\partial_{m_s}f(z^{1/s})={}^{C}D_z^{1/s}(f(z^{1/s}))$ (see~\cite{mi}, Definition 5 and Remark 1).

\subsection{Moment matrix exponential}\label{sec22}
Let $m=(m(p))_{p\ge0}$ be a sequence of positive real numbers. We assume that $m(0)=1$. We define the $m$-exponential of a matrix $M\in\C^{n\times n}$ by the formal expression
\begin{align}\label{exponential}
E_{m}(M)=\sum_{p\geq 0}\frac{M^p}{m(p)}.
\end{align}
If there is no ambiguity on the sequence $m$, we will frequently write $E(M)$ instead of $E_{m}(M)$. Observe that, in principle, the expression of $E(M)$ can be seen as a formal element in $\C[[M]]$, i.e. a formal power series in $M$, or an infinite series of matrices. Note that if $M$ is a diagonal matrix,
$$
M=
\left(\begin{array}{ccccc}
m_{11}&0&0&\cdots&0
\\
0&m_{22}&0&\cdots&0
\\
\vdots&&\ddots&&\vdots
\\
0&\cdots&0&m_{n-1n-1}&0
\\
0&\cdots&0&0&m_{nn}
\end{array}\right),
$$
we have that
$$
M^p=
\left(\begin{array}{ccccc}
m_{11}^p&0&0&\cdots&0
\\
0&m_{22}^p&0&\cdots&0
\\
\vdots&&\ddots&&\vdots
\\
0&\cdots&0&m_{n-1n-1}^p&0
\\
0&\cdots&0&0&m_{nn}^p
\end{array}\right),
$$
therefore
$$
\sum_{p\geq 0}\frac{M^p}{m(p)}=
\left(\begin{array}{ccccc}
\sum_{p\geq 0}\frac{m_{11}^p}{m(p)}&0&0&\cdots&0
\\
0&\sum_{p\geq 0}\frac{m_{22}^p}{m(p)}&0&\cdots&0
\\
\vdots&&\ddots&&\vdots
\\
0&\cdots&0&\sum_{p\geq 0}\frac{m_{n-1n-1}^p}{m(p)}&0
\\
0&\cdots&0&0&\sum_{p\geq 0}\frac{m_{nn}^p}{m(p)}
\end{array}\right),
$$
that is
\begin{align}\label{exp}
E(M)=\left(\begin{array}{ccccc}
E(m_{11})&0&0&\cdots&0
\\
0&E(m_{22})&0&\cdots&0
\\
\vdots&&\ddots&&\vdots
\\
0&\cdots&0&E(m_{n-1n-1})&0
\\
0&\cdots&0&0&E(m_{nn})
\end{array}\right),
\end{align}
at least from a merely formal point of view, and where we have considered in this last step the definition of $E$ for scalars (matrices of dimension 1).

\subsection{Properties of $E(M)$}\label{sec23}

In this section we explore the properties that the operator defined in \eqref{exponential} satisfies. In particular, we compare them with the properties that the classical exponential operator meets when applied to a matrix.

In principle, the definition of the $m$-exponential of a square matrix is given from a formal point of view. In practice, the formal operator $\C^{n\times n}\to \C^{n\times n}[[z]]$ defined by $M\mapsto \sum_{p\ge0}\frac{1}{m(p)}M^pz^p$ will be considered.

\begin{lemma} \label{lemma:properties}
Let $A,B,C\in\C^{n\times n}$. The following formal properties hold:
\begin{enumerate}
\item   If $AB=BA$, then $E(A)E(B)=E(B)E(A)$, when considering the Cauchy product of formal power series.


\item $\partial_mE(Az)=AE(Az)$.

\item If $C$ is invertible and $A=CBC^{-1}$, then $E(A)=CE(B)C^{-1}$, i.e. if the matrices $A$ and $B$ are similar, the same holds for $E(A)$ and $E(B)$.

\item If $O$ is the null matrix, then $m(0)E(O)=I$. In particular, if $m(0)=1$, $E(O)=I$.

\end{enumerate}
\end{lemma}
\begin{proof}
We first prove $(i)$. Note that by the Cauchy product of two formal power series, and since $A$ and B commute one has
\begin{align*}
E(A)E(B)&=\sum_{p\geq 0}\frac{A^p}{m(p)}\sum_{p\geq 0}\frac{B^p}{m(p)}
\\
&= \sum_{p\geq 0}\sum_{n=0}^p\frac{A^n}{m(n)}\frac{B^{p-n}}{m(p-n)}= \sum_{p\geq 0}\sum_{n=0}^p\frac{B^{p-n}}{m(n)}\frac{A^n}{m(p-n)}
\\
&= \sum_{p\geq 0}\sum_{l=0}^p\frac{B^{l}}{m(l)}\frac{A^{p-l}}{m(p-l)}
=\sum_{p\geq 0}\frac{B^p}{m(p)}\sum_{p\geq 0}\frac{A^p}{m(p)}=E(B)E(A).
\end{align*}

 
 The proof of $(ii)$ follows from the definition of the formal moment derivative
 \begin{align*}
\partial_mE(Az)=\partial_m\sum_{p\geq 0}\frac{A^p}{m(p)}z^p=\sum_{p\geq 0}\frac{A^{p+1}}{m(p)}z^p=A\sum_{p\geq 0}\frac{A^p}{m(p)}z^p=AE(Az).
 \end{align*}
 
 As for $(iii)$ it follows directly from the fact that $(CBC^{-1})^p=CB^pC^{-1}$, for every positive integer $p$:
 \begin{align*}
 E(A)=\sum_{p\ge0}\frac{A^p}{m(p)}=\sum_{p\geq 0}\frac{(CBC^{-1})^p}{m(p)}=
 \sum_{p\geq 0}\frac{CB^pC^{-1}}{m(p)}=C\sum_{p\geq 0}\frac{B^p}{m(p)}C^{-1}=CE(B)C^{-1}.
 \end{align*}
 
 Finally, $(iv)$ follows directly from the definition of $E(O)$.
\end{proof}

\begin{remark}
As in the classical setting, it is easy to see that if  $A$  and $B$ do not commute, then the first statement in Lemma  \ref{lemma:properties} does not hold. Besides $E(A+B)\neq E(A)E(B)$ even if  $A$  and $B$ commute, to see this, note that, assuming that $AB=BA$,
\begin{align*}
E(A+B)=\sum_{p\geq 0}\sum_{n=0}^p\frac{p!}{n!(p-n)!m(p)} A^nB^{p-n}
\end{align*}
whereas
\begin{align*}
E(A)E(B)=\sum_{p\geq 0}\sum_{n=0}^p\frac{A^nB^{p-n}}{m(p-n)m(n)}.
\end{align*}
Therefore, in order to be $E(A+B)= E(A)E(B)$ it should hold that
\begin{align*}
\frac{p!}{n!(p-n)!}=\frac{m(p)}{m(n)m(p-n)},\qquad 0\le n\le p,
\end{align*}
which is not true in general. Observe that the previous property implies that the expression of the sequence $m=(m(p))_{p\ge0}$ determines that $m(p)=m(1)^pp!$ by taking $n\ge1$ and $n=p-1$ and from a recursion argument. Therefore, the sequence is equivalent to a Gevrey sequence of order one.  

As a consequence we have  that, in general, we can not deduce that $E(-A)= E(A)^{-1}$, as formal power series with coefficients in the space of $n\times n$ matrices. In fact, in general we have that $E(-A)\neq E(A)^{-1}$. To illustrate  this  consider the following example, let $m(p):=2^{p}$ and $A:=I$, then,
\begin{align*}
E(A)=\sum_{p\geq 0}\frac{I^p}{2^p}=2I,
\end{align*}
whereas 
\begin{align*}
E(-A)=\sum_{p\geq 0}\frac{(-I)^p}{2^p}=I\sum_{p\geq 0}\frac{1}{(-2)^p}=\frac{2}{3}I,
\end{align*}
hence
\begin{align*}
E(A)E(-A)=\frac{4}{3}I.
\end{align*}
\end{remark}

Observe from (iii) in the previous lemma that given $A,B\in\C^{n\times n}$ unitarily similar matrices, then $E(A)$ and $E(B)$ remain unitarily similar matrices. This similarity property will be useful in the sequel to write the solutions of linear systems of moment differential equations in terms of functions of the form $E(\lambda z)$, where $\lambda\in\hbox{spec}(A)$. It is also worth mentioning that the classical identity $\hbox{det}(e^{A})=e^{\hbox{tr}(A)}$ is no longer valid in this general framework unless $m=(B^pp!)_{p\ge0}$ for some $B>0$. This is clear regarding the diagonalizable case together with property (iii) of Lemma~\ref{lemma:properties}. Let $\lambda_1,\ldots,\lambda_n$ be the eigenvalues of $A$ (with possible repetitions). Then, $\hbox{det}(E(A))=\hbox{det}(E(\hbox{diag}(\lambda_1,\ldots,\lambda_n)))=\prod_{j=1}^{n}E(\lambda_j)$, whereas $\hbox{tr}(A)=\hbox{tr}(\hbox{diag}(\lambda_1,\ldots,\lambda_n))=\sum_{j=1}^{n}\lambda_j$ which entails that $E(\hbox{tr}(A))=E(\sum_{j=1}^{n}\lambda_j)$. We finally observe that 
$$\prod_{j=1}^{n}E(\lambda_j)=\sum_{p\ge0}\left(\sum_{j_1+\ldots+j_n=p}\frac{1}{m(j_1)\cdots m(j_n)}\lambda_1^{j_1}\cdots \lambda_n^{j_n}\right)$$
and
$$E(\sum_{j=1}^{n}\lambda_j)=\sum_{p\ge0}\frac{(\lambda_1+\ldots+\lambda_n)^{p}}{m(p)},$$
which coincide if $(\lambda_1+\ldots+\lambda_n)^p=\sum_{j_1+\ldots+j_n=p}\frac{m(p)}{m(j_1)\cdots m(j_n)}\lambda_1^{j_1}\cdots\lambda_n^{j_n}$, for every $p\ge0$. This is only valid if $m=(B^pp!)_{p\ge0}$ for some $B>0$. The general case can be considered in a similar manner taking into account the Jordan decomposition of $A$ and the properties of the action of the operator $E(\cdot)$ on Jordan blocks, described in Section~\ref{secjordan}.

Regarding convergence properties, we consider a normalized submultiplicative norm $\left\|\cdot\right\|$ in $\C^{n\times n}$. We will be interested in the convergence of the formal power series $E(Az)\in\C^{n\times n}[[z]]$.

It is natural for a sequence $m=(m(p))_{p\ge0}$ to satisfy the property
\begin{equation}\label{e433}
\liminf_{p\to\infty}m(p)^{1/p}=+\infty
\end{equation}
when working with ultraholomorphic and/or ultradifferentiable functions, i.e. spaces of functions whose derivatives are subjected to bounds given in terms of $m$. In this respect, we refer to the definition of weight sequence in Section 3~\cite{jss1}. We also refer to~\cite{balser}, Appendix B, for a review on holomorphic functions with values in complex Banach spaces and their main properties.

\begin{lemma}\label{lema450}
Let $m=(m(p))_{p\ge0}$ be a sequence of positive real numbers such that (\ref{e433}) holds. Then, for every $A\in\C^{n\times n}$ the formal power series $E(Az)$ defines an entire funtion (with values in the Banach space $\C^{n\times n}$). In particular, $E(A)$ converges absolutely. In addition to this, for every $z\in\C$ and every normalized submultiplicative norm $\left\|\cdot\right\|$ on $\C^{n\times n}$ one has that 
$$\left\|E(Az)\right\|\le E(\left\|A\right\||z|),\qquad z\in\C$$
\end{lemma}
\begin{proof}
It is straight to check that 
\begin{equation}\label{e451}
\left\|E(Az)\right\|\le \sum_{p\ge0}\frac{1}{m(p)}\left\|A\right\|^{p}|z|^{p}.
\end{equation}
We observe that the radius of convergence of the series on the right-hand side of the previous inequality is $\rho$, where
$$\frac{1}{\rho}=\left\|A\right\|\limsup_{p\to\infty}\frac{1}{m(p)^{1/p}}=\left\|A\right\|\left(\liminf_{p\to\infty}(m(p))^{1/p}\right)^{-1}=0.$$
For the second part of the proof, let $z\in\C$. We observe that the right-hand side of (\ref{e451}) coincides with $E(\left\|A\right\||z|)$.
\end{proof}

The construction of $E(A)^{-1}$ mentioned above is determined recursively as follows.

\begin{lemma}
Let $A\in\C^{n\times n}$ and $m=(m(p))_{p\ge0}$ be a sequence of positive real numbers such that (\ref{e433}) holds. Then, $E(Az)$ is an invertible matrix and 
$$E(Az)^{-1}=\sum_{p\ge0}\frac{\phi_p}{m(p)}A^pz^p\in\C^{n\times n}\{z\},$$
with $\phi_0=1$, and $\phi_{p}=-\sum_{j=0}^{p-1}\frac{m(p)}{m(j)m(p-j)}\phi_j$ for every $p\ge1$.
\end{lemma}
\begin{proof}
Write a generic formal power series $B(z)=\sum_{p\ge0}\frac{1}{m(p)}B_pz^p\in\C^{n\times n}[[z]]$. If $B(z)E(Az)=I$, then one can prove inductively that $B_p=\phi_p A^p$, for some $\phi_p\in\C$ which satisfies that $\phi_0=1$ and $\sum_{j=0}^{p}\frac{1}{m(j)m(p-j)}\phi_j=0$ for all $p\ge1$. The formal existence of $E(Az)^{-1}$ is then guaranteed. Indeed, as $E(Az)$ defines an entire function (with values in the Banach space $\C^{n\times n}$) and $E(Az)|_{z=0}=I$, there exists a neighborhood of the origin in which $E(Az)^{-1}$ defines an analytic function. More precisely, the radius of convergence is determined by the modulus of the first zero of the function $z\mapsto E(Az)$. An alternative inductive argument can be followed to check the existence large enough $c,b>0$ such that $\phi_p/m(p)\le c b^p$ for every $p\ge0$. 
\end{proof}

Observe that in the case that $m=(p!)_{p\ge0}$, then $\sum_{j=0}^{p}\frac{m(p)}{m(j)m(p-j)}\phi_j=0$ for all $p\ge1$, with $\phi_j=(-1)^j$ for every $j\ge0$ and therefore $E(Az)^{-1}=\exp(Az)^{-1}=\exp(-Az)=E(-Az)$. 

The property (ii) stated in Lemma~\ref{lemma:properties} can be reconsidered in the convergent setting. The next result is a direct consequence of Lemma~\ref{lema450}.

\begin{lemma}
Let $A\in\C^{n\times n}$, and $m$ be a sequence of positive real numbers such that (\ref{e433}) holds. Then, $\partial_{m}E(Az)$ determines an entire function such that $\left\|\partial_mE(Az)\right\|\le \left\|A\right\|E(\left\|A\right\| |z|)$, for every $z\in\C$.
\end{lemma}

\subsection{Jordan blocks}\label{secjordan}
A Jordan block, $J_{\lambda,i}$ is a $i\times i$ upper triangular matrix such that it has the same number, $\lambda\in \mathbb{C}$, in the diagonal an every entry on the superdiagonal is 1, that is
$$
J_{\lambda,i}=
\left(\begin{array}{ccccc}
\lambda&1&0&\cdots&0
\\
0&\lambda&1&\ddots&\vdots
\\
0&\ddots&\ddots&\ddots&0
\\
\vdots&\ddots&0&\lambda&1
\\
0&\cdots&0&0&\lambda
\end{array}
\right)=\lambda I_i+N_i,
$$
where $I_i$ is the  $i\times i$ identity matrix and 
$$
N_i=\left(\begin{array}{ccccc}
0&1&0&\cdots&0
\\
0&0&1&\ddots&\vdots
\\
0&\ddots&\ddots&\ddots&0
\\
\vdots&\ddots&0&0&1
\\
0&\cdots&0&0&0
\end{array}
\right)
$$
is a nilpotent matrix such that $N_i^p$ is the null matrix for all $p\geq i$. Therefore,
\begin{align*}
J_{\lambda,i}^p=(\lambda I_i+N_i)^p=\sum_{h=0}^{\min\{i-1,p\}}\binom{p}{h}\lambda^{p-h}I_iN_i^h=\sum_{h=0}^{\min\{i-1,p\}}\binom{p}{h}\lambda^{p-h}N_i^h.
\end{align*}
Hence,
\begin{align}\label{expJ}
E(J_{\lambda,i})&=\sum_{p\geq 0}\frac{J_{\lambda,i}^p}{m(p)}
\\
\nonumber
&=
\sum_{p\geq 0}\sum_{h=0}^{\min\{i-1,p\}}\binom{p}{h}\frac{\lambda^{p-h}}{m(p)}N_i^h=\sum_{h= 0}^{i-1}N_i^h\sum_{p\geq h}\binom{p}{h}\frac{\lambda^{p-h}}{m(p)}
\\
\nonumber
&=
\left(\begin{array}{ccccc}
\Delta_{0}E(\lambda,1)&\Delta_{1}E(\lambda,1)&\cdots&\Delta_{i-2}E(\lambda,1)&\Delta_{i-1}E(\lambda,1)
\\
0&\Delta_{0}E(\lambda,1)&\Delta_{1}E(\lambda,1)&\cdots&\Delta_{i-2}E(\lambda,1)
\\
\\
\vdots&\ddots&\ddots&\ddots&\vdots
\\
\\
0&\cdots&0&\Delta_{0}E(\lambda,1)&\Delta_{1}E(\lambda,1)
\\
0&\cdots&0&0&\Delta_{0}E(\lambda,1)
\end{array}\right).
\end{align}
 
Here, we have used the following definition.

\begin{definition}[Definition 5,~\cite{Lastra22}]\label{defidelta}
Let $\lambda\in\C$ and $h\ge0$. We define
$$\Delta_{h}E(\lambda,z)=\sum_{p\ge h}\binom{p}{h}\frac{\lambda^{p-h}z^p}{m(p)}.$$
\end{definition}
\section{Main results}
In this section, we provide the general solution to (\ref{problem}) in terms of the generalized exponential matrix, in the spirit of the classical exponential matrix. We distinguish two cases: 1) the matrix of the system (\ref{problem}) is diagonalizable; 2) the martix is not diagonalizable. 

In this section, we assume $m=(m(p))_{p\ge0}$ is a sequence of moments constructed from a kernel of generalized summability. More precisely, let $\mathbb{M}=(M_p)_{p\ge0}$ be a strongly regular sequence admitting a nonzero proximate order. For a pair of associated kernels for generalized summability $e$ and $E$, the sequence $m$ is determined by (\ref{e336}). Given a (lc) sequence $\mathbb{M}$ (in particular a strongly regular sequence), it is straight to check from a convexity argument that the sequence $(M_p^{1/p})_{p\ge0}$ is monotone increasing. The fact that $\sup_{p\ge0}M_p^{1/p}<\infty$ is equivalent to the corresponding class of ultradifferentiable functions to be contained in the space of analytic functions~\cite{cm}, which leads to trivial situations. For this reason, we assume that $\lim_{p\to\infty}M_p^{1/p}=+\infty$. Regarding~\cite{sa0}, Proposition 5.8, one can guarantee the existence of $A,B>0$ such that 
$$A^pM_p\le m(p)\le B^pm_p,\quad p\ge0,$$
which entails that (\ref{e433}) holds. We will also assume this condition holds from now on. Observe this is the case of any Gevrey sequence in particular.

We also recall the following results from~\cite{Lastra22} for the sake of completeness.

\begin{lemma}[Lemma 2,\cite{Lastra22}]
The set of solutions to (\ref{problem}) is a subspace of $(\mathcal{O}(\C)^{n}$ of dimension $n$.
\end{lemma}

\begin{lemma}[Lemma 5,~\cite{Lastra22}]
For all $\lambda\in\C$ and $h\ge0$, $\Delta_hE(\lambda,z)\in\mathcal{O}(\C)$, with $\Delta_hE(\lambda,z)$ defined in Definition~\ref{defidelta}.
\end{lemma}

\subsection{A diagonalizable}
If $A$ is a diagonalizable matrix there exist $D$, diagonal matrix, and $P$ an invertible matrix such that
$$
A=PDP^{-1}.
$$
Then,
$$
D=\left(\begin{array}{ccccc}
\lambda_{1}&0&0&\cdots&0
\\
0&\lambda_{2}&0&\cdots&0
\\
\vdots&&\ddots&&\vdots
\\
0&\cdots&0&\lambda_{n-1}&0
\\
0&\cdots&0&0&\lambda_{n}
\end{array}\right)
\quad
\textrm{ and }
\quad
P=\left(\begin{array}{ccc}
v_{1}&\cdots&v_n
\end{array}\right),
$$
where $\lambda_i$ is an eigenvalue with associated eigenvector $v_i$ (column vector), $i\in\{1,\ldots,n\}$, note that the $\lambda_i$ are not necessarily different. 

Since $E(A)=\sum_{p\geq 0}\frac{A^p}{m(p)}$,
in view of \eqref{exp} we have that
$$
E(A)=PE(D)P^{-1}=P\left(\begin{array}{ccccc}
E(\lambda_{1})&0&0&\cdots&0
\\
0&E(\lambda_{2})&0&\cdots&0
\\
\vdots&&\ddots&&\vdots
\\
0&\cdots&0&E(\lambda_{n-1})&0
\\
0&\cdots&0&0&E(\lambda_{n})
\end{array}\right)P^{-1}.
$$

\begin{theorem}\label{teo1}
Let $A\in \mathbb{C}^{n\times n}$ be a diagonalizable matrix,  the general solution of \eqref{problem} is given by
\begin{align}\label{solution}
y(z)=E(Az)v^c,
\end{align}
where  $v^c$ is an n-dimensional constant column vector. 

\end{theorem}
\begin{proof}
 Let $\{\lambda_i\}_{1\leq i\leq n}$ be the set of eigenvalues of $A$ with associated eigenvectors $\{v_i\}_{1\leq i\leq n}$. Then,
expanding \eqref{solution}, we have
\begin{align*}
E(Az)v^c=\left(\begin{array}{ccc}
E(\lambda_{1}z)v_{1}&\cdots&E(\lambda_{n}z)v_n
\end{array}\right)\left(\begin{array}{ccc}
\widetilde{v}_{1}&\cdots&\widetilde{v}_n
\end{array}\right)\left(\begin{array}{c}
v_{1}^c\\
\vdots
\\
v_n^c
\end{array}\right),
\end{align*}
where $\left(\begin{array}{ccc}
\widetilde{v}_{1}&\cdots&\widetilde{v}_n
\end{array}\right)$ is  a vector column matrix such that 
$$
\left(\begin{array}{ccc}
\widetilde{v}_{1}&\cdots&\widetilde{v}_n
\end{array}\right)^{-1}=\left(\begin{array}{ccc}
{v}_{1}&\cdots&{v}_n
\end{array}\right).
$$
Note that
\begin{align*}
\left(\begin{array}{ccc}
\widetilde{v}_{1}&\cdots&\widetilde{v}_n
\end{array}\right)\left(\begin{array}{c}
v_{1}^c\\
\vdots
\\
v_n^c
\end{array}\right)=
\widetilde{v}_{1}v_1^c+\cdots+\widetilde{v}_n
v_n^c,
\end{align*}
where remember that $\widetilde{v}_{1},\,\cdots,\,\widetilde{v}_n$ are $n$ linear independent column vectors.

Then, if we define 
$$
\left(\begin{array}{c}
c_1\\
\vdots
\\
c_n
\end{array}\right):=
\widetilde{v}_{1}v_1^c+\cdots+\widetilde{v}_n
v_n^c,
$$
we have that
\begin{align*}
E(Az)v^c=\left(\begin{array}{ccc}
E(\lambda_{1}z)v_{1}&\cdots&E(\lambda_{n}z)v_n
\end{array}\right)\left(\begin{array}{c}
c_1\\
\vdots
\\
c_n
\end{array}\right)=\sum_{i=1}^nc_iE(\lambda_iz)v_i.
\end{align*}
Besides, note that for any choice of $v^c$ the vector $\left(\begin{array}{c}
c_1\\
\vdots
\\
c_n
\end{array}\right)$ is uniquely determined since 
$$
\textrm{det}\left(\begin{array}{ccc}
\widetilde{v}_{1}&\cdots&\widetilde{v}_n
\end{array}\right)\neq 0.
$$
Therefore, it is sufficient to apply \cite[Theorem 1]{Lastra22} to conclude the proof.

\end{proof}

\subsection{A not diagonalizable}
Assuming that $A$ is a non-diagonalizable matrix 
there exist $\lambda_i$, $1\leq i\leq k$, eigenvalues of $A$ with algebraic multiplicity $m_i$, $1\leq i\leq k$ ($\sum_{i=1}^km_i=n$), consider $\{v_{i,j}\}_{1\leq j\leq m_i}$ the associated generalized eigenvectors of each $\lambda_i$. Then, there exits a Jordan matrix $J$, and an invertible matrix $P$ whose columns are determined by generalized eigenvectors associated to $\lambda_i$ such that
$$
A=PJP^{-1},
$$

$$
J=\left(\begin{array}{cccccccc}
J_{\lambda_1,l_{1,1}}&0&0&\cdots&&&0
\\
0&\ddots&&&&&\vdots
\\
\vdots&&J_{\lambda_1,l_{1,q_1}}&0&\cdots&&0
\\
0&\cdots&0&\ddots&&&0
\\
\vdots&&&&\ddots&&\vdots
\\
0&\cdots&&&0&J_{\lambda_{k},l_{k,1}}&
\\
\vdots&&&&&\,\,\quad\ddots&
\\
0&\cdots&&&&0&J_{\lambda_k,l_{k,q_k}}
\end{array}\right)
$$
where each $J_{\lambda_i,l_{i,s}}$, $1\leq i\leq k$, $1\leq s\leq q_i$ ($\sum_{s=1}^{q_i}l_{i,s}=m_i$), is a Jordan block of size $l_{i,s}$.
\begin{theorem}\label{teo2}
Let $A\in \mathbb{C}^{n\times n}$ be a non-diagonalizable matrix. The general solution of \eqref{problem} is given by
\begin{align}\label{solutionnond}
y(z)=E(Az)v^c,
\end{align}
where  $v^c$ is an n-dimensional constant column vector. 

\end{theorem}
\begin{proof}
Let $\lambda_{i}$, ${1\leq i\leq k}$ be the eigenvalues of $A$, with algebraic multiplicity $m_i$,  $\sum_{i=1}^km_i=n$. For each  eigenvalue $\lambda_i$ we consider its associated eigenvectors $\{v_{i,l_{i,s}}^j\}_{1\leq s\leq q_i}$, $1\leq j\leq l_{i,s}$.

Expanding $E(Az)$, we have
\begin{align*}
&E(Az)
\\
&=\left(\begin{array}{ccccccc}
\Delta_0E(\lambda_{1},z)v_{1,l_{1,1}}^1&\!\!\!\cdots\!\!\!&\sum_{h=1}^{l_{1,1}}\Delta_{l_{1,1}-h}E(\lambda_{1},z)v_{1,l_{1,1}}^{l_{1,1}}&\!\!\!\cdots\!\!\!&\Delta_0E(\lambda_{k},z)v_{k,l_{k,q_k}}^1&\!\!\!\cdots\!\!\!&\sum_{h=1}^{l_{k,q_k}}\Delta_{l_{k,q_k}-h}E(\lambda_{k},z)v_{k,l_{k,q_k}}^{l_{k,q_k}}
\end{array}\right)P^{-1},
\end{align*}
where 
$$
P^{-1}=\left(\begin{array}{ccccccccc}
{v}_{1,l_{1,1}}^1&\cdots&{v}_{1,l_{1,1}}^{l_{1,1}}&\cdots&{v}_{k,l_{k,q_k}}^1\cdots&{v}_{k,l{k,q_k}}^{l_{k,q_k}}
\end{array}\right)^{-1}=:\left(\begin{array}{ccc}
\widetilde{v}_{1}&\cdots&\widetilde{v}_n
\end{array}\right),
$$
where  $\left(\begin{array}{ccc}
\widetilde{v}_{1}&\cdots&\widetilde{v}_n
\end{array}\right)$ is a vector column matrix.

Besides
\begin{align*}
P^{-1}\left(\begin{array}{c}
v_{1}^c\\
\vdots
\\
v_n^c
\end{array}\right)=
\widetilde{v}_{1}v_1^c+\cdots+\widetilde{v}_n
v_n^c,
\end{align*}
where remember that $\widetilde{v}_{1},\,\cdots,\,\widetilde{v}_n$ are $n$ linear independent column vectors.
Then, for any choice of $v^c$, the vector 
$$
\left(\begin{array}{c}
c_{1,l_{1,1}}^1\\
\vdots
\\
c_{1,l_{1,1}}^{l_{1,1}}
\\
\vdots
\\
c_{k,l_{k,q_k}}^1
\\
\vdots
\\
c_{k,l_{k,q_k}}^{l_{k,q_k}}
\end{array}\right)
:=
\widetilde{v}_{1}v_1^c+\cdots+\widetilde{v}_n
v_n^c,
$$
is uniquely determined since 
$$
\textrm{det}\left(\begin{array}{ccc}
\widetilde{v}_{1}&\cdots&\widetilde{v}_n
\end{array}\right)\neq 0.
$$
Besides, recall that $\sum_{i=1}^k\sum_{s=1}^{q_i}l_{i,s}=\sum_{i=1}^km_i=n$.

Gathering the above estimates, we have that $y(z)=E(Az)v^c$ coincides with
\begin{multline*}\left(\begin{array}{cccc}
\Delta_0E(\lambda_{1},z)v_{1,l_{1,1}}^1&\!\!\!\cdots\!\!\!&\sum_{h=1}^{l_{1,1}}\Delta_{l_{1,1}-h}E(\lambda_{1},z)v_{1,l_{1,1}}^{l_{1,1}}&\!\!\!\cdots\end{array}\right.\\
\left.\begin{array}{cccc}\!\!\!\cdots\!\!\! &\Delta_0E(\lambda_{k},z)v_{k,l_{k,q_k}}^1&\!\!\!\cdots\!\!\!&\sum_{h=1}^{l_{k,q_k}}\Delta_{l_{k,q_k}-h}E(\lambda_{k},z)v_{k,l_{k,q_k}}^{l_{k,q_k}}
\end{array}\right)
\left(\begin{array}{c}
c_{1,l_{1,1}}^1\\
\vdots
\\
c_{1,l_{1,1}}^{l_{1,1}}
\\
\vdots
\\
c_{k,l_{k,q_k}}^1
\\
\vdots
\\
c_{k,l_{k,q_k}}^{l_{k,q_k}}
\end{array}\right)
\\
=\sum_{i=1}^k\sum_{s=1}^{q_i}\sum_{j=1}^{l_{i,s}}c_{i,l_{i,s}}^j\sum_{h=1}^j\Delta_{j-h}E(\lambda_i,z)v_{i,l_{i,s}}^j.
\end{multline*}

Therefore, applying \cite[Theorem 2]{Lastra22} we conclude the proof.
\end{proof}

\begin{definition}
Let $A\in\C^{n\times n}$ and consider the linear system of moment differential equations (\ref{problem}). A matrix $X=X(z)\in(\mathcal{O}(\C))^{n\times n}$ is said to be a fundamental matrix associated to (\ref{problem}) if the columns of $X$ determine a set of linearly independent solutions of (\ref{problem}).
\end{definition}

As a consequence of Theorem~\ref{teo1} and Theorem~\ref{teo2} we derive the following result.

\begin{corollary}
Let $A\in\C^{n\times n}$ and consider the linear system of moment differential equations (\ref{problem}). Then, the matrix $E(Az)$ is a fundamental matrix associated to (\ref{problem}). 
\end{corollary}

\begin{corollary}
Let $X(z)\in(\mathcal{O}(\C))^{n\times n}$ be a fundamental matrix associated to (\ref{problem}). Then, it holds that
$$E(Az)=X(z)X(0)^{-1}.$$
\end{corollary}
\begin{proof}
Given two fundamental solutions $X(z),Y(z)\in(\mathcal{O}(\C))^{n\times n}$ of (\ref{problem}), it is straight the existence of an invertible matrix $C\in\C^{n\times n}$ such that $Y(z)=X(z)C$. This is a consequence of the columns of $X$ and the columns of $Y$ being basis of the vector space of solutions to (\ref{problem}), so $C$ represents the change of coordinates matrix. This entails that $E(Az)=X(z)C$. The evaluation at $z=0$ yields $E(0)=I=X(0)C$. Therefore, $C=X(0)^{-1}$. Recall that $X(0)$ is an invertible matrix which represents, in columns a set of $n$ linearly independent initial conditions at the origin of (\ref{problem}), given by a basis of $\C^n$.
\end{proof}

As in the classical situation, the computation of $E(Az)$ is more easy to handle than the linear combination (see the proof of Theorem 2~\cite{Lastra22}) , providing the whole information in a compact expression.  

\begin{example}
Let $A\in\C^{n\times n}$ be given by
$$A=\left(\begin{array}{ccc}
1&0&1\\
1&2&0\\
0&0&1
\end{array}\right),$$
and let $m$ be a sequence of positive real numbers which satisfies (\ref{e433}). We consider the linear system of moment differential equations $\partial_my=Ay$. It is straight to check that the characteristic polynomial associated with $A$ has $r_1=1$ as a double root, and $r_2=2$ as a single root. The Jordan canonical form associated to $A$ is given by $J$ and matrix $P$
$$J=\left(\begin{array}{ccc}
1&1&0\\
0&1&0\\
0&0&2
\end{array}\right),\qquad P=\left(\begin{array}{ccc}
1&-1&0\\
-1&0&1\\
0&1&0
\end{array}\right),$$
with $P^-1AP=J$. The general solution of $\partial_my=Ay$ is $y(z)=E(Az)$. On one hand, one can prove by induction that
$$A^p=\left(\begin{array}{ccc}
1&0&p\\
2^p-1&2^p&2^p-p-1\\
0&0&1
\end{array}\right),$$
which entails that $E(Az)=\sum_{p\ge0}\frac{1}{m(p)}A^pz^p$ coincides with
$$\left(\begin{array}{ccc}
\sum_{p\ge0}\frac{1}{m(p)}z^p&0&\sum_{p\ge0}\frac{p}{m(p)}z^p\\
\sum_{p\ge0}\frac{1}{m(p)}(2z)^p-\sum_{p\ge0}\frac{1}{m(p)}z^p& \sum_{p\ge0}\frac{1}{m(p)}(2z)^p& \sum_{p\ge0}\frac{1}{m(p)}(2z)^p-\sum_{p\ge0}\frac{p}{m(p)}z^p-\sum_{p\ge0}\frac{1}{m(p)}z^p\\
0&0& \sum_{p\ge0}\frac{1}{m(p)}z^p
\end{array}\right).$$
On the other hand, one obtains that $E(Az)=PE(J_1(z))P^{-1}$, with $J_1(z)\in(\mathcal{O}(\C))^{3\times 3}$ being the blocked matrix
$$E(J_1(z))=\left(\begin{array}{cc}
J_{11}(z)&0\\
0&J_{22}
\end{array}\right),$$
with 
$$J_{11}(z)=\left(\begin{array}{cc}
\Delta_0E(1,1)&\Delta_1E(1,1)\\
0&\Delta_0E(1,1)
\end{array}\right)=\left(\begin{array}{cc}
E(z)&\sum_{p\ge0}\frac{p}{m(p)}z^p\\
0&E(z)
\end{array}\right),\qquad J_{22}(z)=(\Delta_0E(2,1))=(E(2z)).$$
One can directly check that $PE(J_1(z))P^{-1}$ coincides with $E(Az)$ computed above.
\end{example}

\begin{remark}
The particular case of $m_k=(\Gamma(1+\frac{p}{k})_{p\ge0}$ for some fixed positive integer $k$ has  been previously studied. For the initial case $k=1$, the moment derivation is reduced to the classical derivation. If $k\ge2$, the moment derivative is quite related to Caputo fractional derivative, as stated above, with $e(z)=kz^{k}\exp(-z^k)$, and $E(z)$ being Mittag-Leffler function $E_{1/k}(z)=\sum_{p\ge0}\frac{1}{\Gamma(1+p/k)}z^p$. In~\cite{boritr}, the authors provide the solution to this problem in the Riemann-Liouville and Caputo derivative, and describe the solution in terms of a fractional exponential matrix. In this sense, our result particularizes to Theorem 4 in~\cite{boritr} when chosing $m_k$ to be the sequence of moments.
\end{remark}

\subsection{Beyond moderate growth}

In this last paragraphs, we aim to say some words about the application of the previous results to the $q-$Gevrey setting, i.e. when the sequence $m$ is of the form $m_q=([p]_q!)_{p\ge0}$. We recall that, given $q\in\R$ with $q\neq0,1$, the $q-$factorial numbers are defined by $[0]_q!=1$ and $[p]_q!=[p]_q\cdot [p-1]_q\cdot\ldots\cdot[1]_q$, for every positive integer $p$, where $[k]_q$ stands for the $q-$number $[k]_q=1+q+\ldots+q^{k-1}$.

This sequence is no longer a strongly regular sequence due to (mg) condition is not satisfied, so the previous results are no longer available. However, the moment derivative associated to this sequence is of great importance due to $q-$derivative coincides with such moment derivation. Indeed,
$$\partial_{m_q}f(z)=D_qf(z)=\frac{f(qz)-f(z)}{qz-z},$$
which converges to the usual derivative when $q\to 1$.

Let $q>1$. The sequence $m_q$ is quite related to the sequence $(q^{p(p-1)/2})_{p\ge0}$ as comparable sequences. The latter being the sequence of moments (in the sense of (\ref{e336}) associated to the kernel $e(z)=\sqrt{2\pi\ln(q)}\exp\left(\frac{\ln^2(\sqrt{q}z)}{2\ln(q)}\right)$, and also to the kernel $e(z)=\ln(q)\Theta_{1/q}(z)$, where $\Theta_{1/q}$ stands for Jacobi Theta function $\Theta_{1/q}(z)=\sum_{p\in\mathbb{Z}}\frac{1}{q^{\frac{p(p-1)}{2}}}z^p$, which is holomorphic on $\C\setminus\{0\}$ with an essential singularity at the origin. Indeed, observe that this sequence of moments does not satisfy Carleman condition for determinacy of the Stieltjes moment problem (see~\cite{ca1}). 

Following~\cite{ra}, let us consider the function $\exp_{q}(z)=\sum_{p\ge0}\frac{1}{[p]_q!}z^p$, which turns out to be an entire function. The sequence $m_q$ satisfies property (\ref{e433}), and the formal and analytic results obtained in Section~\ref{sec22} and Section~\ref{sec23} hold for $E(z)=\exp_q(z)$. In addition to this, the construction of the solutions to
$$D_qy=Ay$$ 
of Theorem 2 and Theorem 3~\cite{Lastra22} may apply. Therefore, $y(z)=\exp_{q}(Az)$ is a fundamental matrix associated to the previous problem.

\begin{example}
Let $A\in\C^{n\times n}$ be given by
$$A=\left(\begin{array}{ccc}
0&1&1\\
-1&2&1\\
1&-1&1
\end{array}\right),$$
We consider the linear system of moment differential equations $D_qy=Ay$. It is straight to check that the characteristic polynomial associated with $A$ has $r_1=1$ as a triple root. The Jordan canonical form associated to $A$ is given by $J$ and matrix $P$
$$J=\left(\begin{array}{ccc}
1&1&0\\
0&1&1\\
0&0&1
\end{array}\right),\qquad P=\left(\begin{array}{ccc}
1&0&1\\
1&0&0\\
0&1&1
\end{array}\right),$$
with $P^-1AP=J$. The general solution of $D_qy=Ay$ is $y(z)=exp_q(Az)$. One has that
$$A^p=\frac{1}{2}\left(\begin{array}{ccc}
p^2-3p+2&-(p-3)p&2p\\
(p-3)p&-p^2+3p+2&2p\\
2p&-2p&2
\end{array}\right),$$
which entails that $\exp_q(Az)=\sum_{p\ge0}\frac{1}{[p]_{q}!}A^pz^p$ is given by
$$\left(\begin{array}{ccc}
\sum_{p\ge0}\frac{1}{[p]_q!}\frac{1}{2}(p^2-3p+2)z^p & \sum_{p\ge0}\frac{1}{[p]_q!}\frac{-1}{2}(p-3)pz^p & \sum_{p\ge0}\frac{1}{[p]_q!}pz^p
\\
\sum_{p\ge0}\frac{1}{[p]_q!}\frac{1}{2}(p-3)pz^p & \sum_{p\ge0}\frac{1}{[p]_q!}\frac{1}{2}(-p^2+3p+2)z^p & \sum_{p\ge0}\frac{1}{[p]_q!}pz^p
\\
\sum_{p\ge0}\frac{1}{[p]_q!}pz^p & \sum_{p\ge0}\frac{-1}{[p]_q!}pz^p & \sum_{p\ge0}\frac{1}{[p]_q!}z^p
\end{array}\right).$$

Observe that $P(\exp_q(Jz))P^{-1}$ coincides with $\exp_q(Az)$,
$$\exp_q(Jz)=\left(\begin{array}{ccc}
\Delta_0E(1,1) & \Delta_1E(1,1) & \Delta_2E(1,1)\\
0 & \Delta_0E(1,1) & \Delta_1E(1,1)\\
0 & 0 & \Delta_0E(1,1)
\end{array}\right)=\left(\begin{array}{ccc}
\exp_q(z) & \sum_{p\ge0}\frac{p}{[p]_q!}z^p & \sum_{p\ge0}\binom{p}{2}\frac{1}{[p]_q!}z^p\\
0 & \exp_q(z) & \sum_{p\ge0}\frac{p}{[p]_q!}z^p\\
0 & 0 & \exp_q(z) 
\end{array}\right).$$
\end{example}


\end{document}